\documentclass[12pt, a4paper]{article}
\usepackage[pagewise]{lineno}
\usepackage{graphicx}
\usepackage{amssymb}
\usepackage{latexsym, bm}
\usepackage{multicol}
\usepackage{indentfirst}
\usepackage{amssymb,amsfonts}
\usepackage{amsmath}
\usepackage{cases}
\usepackage[noadjust]{cite}
\usepackage[colorlinks,citecolor=red]{hyperref}
\newtheorem{theorem}{Theorem}
\newtheorem{lemma}{Lemma}

\newtheorem{definition}{Definition}

\newtheorem{remark}{Remark}

\usepackage{amssymb}
\numberwithin{equation}{section}
\numberwithin{theorem}{section}
\numberwithin{remark}{section}
\numberwithin{definition}{section}
\numberwithin{lemma}{section}
\numberwithin{corollary}{section}
\numberwithin{proposition}{section}
\title{Energy conservation for the non-resistive MHD equations with physical boundaries}
\author{Wenke Tan\footnote{tanwenkeybfq@163.com}\quad Fan Wu\footnote{wufan0319@yeah.net}\\
{\small Key Laboratory of Computing and Stochastic Mathematics (Ministry of Education),}\\
{\small School of Mathematics and Statistics, Hunan Normal University,}\\
{\small Changsha, Hunan 410081, China}\\
}

\date{}

\begin{document}
\maketitle
{\bf Abstract:}
In this paper, we study the energy equality for weak solutions to the non-resistive MHD equations with physical boundaries. Although the equations of magnetic field $b$ are of hyperbolic type, and the boundary effects are considered, we still prove the global
energy equality provided that
$
	u \in L^{q}_{loc}\left(0, T ; L^{p}(\Omega)\right)
	\text { for any } \frac{1}{q}+\frac{1}{p} \leq \frac{1}{2}, \text { with } p \geq 4,\text{ and }  b \in L^{r}_{loc}\left(0, T ; L^{s}(\Omega)\right) \text { for any } \frac{1}{r}+\frac{1}{s} \leq \frac{1}{2}, \text { with } s \geq 4  $.
In particular, compared with the existed
results, we do not require any boundary layer assumptions and additional conditions on
the pressure $P$. Our result requires the regularity of boundary $\partial\Omega$ is only Lipschitz which is the minimum requirement to make the boundary condition $b\cdot n$ sense.
The proof is based on the important properties of weak solutions of the nonstationary
Stokes system and the separate mollification of weak solutions from the boundary effect by considering a
non-standard local energy equality and transform the boundary effects into the estimates
of the gradient of cut-off functions.
% Then, by establishing a sharp $L^2L^2$ estimate for pressure
%$P$, we use zero boundary conditions of $u$ to inhibit the boundary effect and obtain global
%energy equality by choosing suitable cut-off functions.

\medskip
{\bf Mathematics Subject Classification (2010):} \  76W05, 76B03, 35Q35.
\medskip

{\bf Keywords:}  Non-resistive MHD equations; Weak solutions; Energy conservation;  Physical boundaries
\section{Introduction}
In this paper, we are concerned with the energy conservation of the weak solutions to the incompressible non-resistive MHD equations with zero external force in $\Omega\subseteq\mathbb{R}^d$ $(d\geq 2)$:
\begin{equation}\label{1.1}
\left\{
             \begin{array}{lr}
             \partial_t u+(u\cdot\nabla)u-\mu\Delta u+\nabla P=(b\cdot\nabla)b,& \\
            \partial_t b+(u\cdot\nabla)b=(b\cdot\nabla)u,&\\
             \nabla\cdot u=\nabla\cdot b=0,&\\
             u|_{\partial\Omega}=0,b\cdot n|_{\partial\Omega}=0,\\
             u(x,0)=u_{0}(x),b(x,0)=b_{0}(x),&
\end{array}
\right.
\end{equation}
where the spatial variable $x \in \Omega \subset \mathbb{R}^{d}$ with $\Omega$ being a bounded
 Lipshitz domain with boundary $\partial \Omega$ and $n$ being the outward unit normal vector filed to the boundary, and $u=u(x, t)$ is the fluid velocity field, $b=b(x, t)$ is the magnetic field, the positive constant
$\mu$ is the viscosity coefficient, $P$ is the scalar pressure, and $u_{0},b_{0}$ are the prescribed initial data satisfying the compatibility condition $\nabla\cdot u_{0}=\nabla\cdot b_{0}$ = 0 in the distributional sense. Physically, equations \eqref{1.1} govern the dynamics of the velocity and magnetic fields in electrically conducting fluids, such as plasmas, liquid metals, and salt water. More details can be found in \cite{CH,TEMAM}.

System \eqref{1.1} can be applied to model plasmas when the plasmas are strongly collisional or when the resistivity due to these collisions is extremely small. It often applies to the case when one is interested in the k-length scales that are much longer than the ion skin depth and the Larmor radius perpendicular to the field, long enough along the field to ignore the Landau damping, and time scales much longer than the ion gyration time \cite{CH,LZ}. In this case, the study of well-posedness will become more difficult. Even the global existence of Leray-Hopf weak solutions is not known due to the lack of suitable strong convergence in magnetic field $b$. In spite of the difficulties due to the lack of
magnetic diffusion, significant progress has been made in \cite{LXZ,ZZF},  under the assumption that the initial data $(u_0, b_0)$ is close to the equilibrium state $\left(0, (1,0)^T\right)$, the global well-posedness was proven. For more progress in this direction, one may refer to the survey paper \cite{LINFANGHUA} and the references therein. In addition, for 2D non-resistive MHD equations \eqref{1.1}, Jiu-Niu \cite{JIUNIU} first proved the local existence and uniqueness of the classical solution in the Sobolev space $H^s$ with $s\geq3$. Fefferman-McCormick-Robinson-Rodrigo \cite{FEFF} were able to weaken the regularity assumption to $(u_0,b_0)\in H^s$ with $s>\frac{d}{2}$ and obtained local-in-time existence of strong solutions to \eqref{1.1} in $\mathbb{R}^d$, $d=2,3$. And then, they made a further improvement by assuming $\left(u_{0}, b_{0}\right) \in H^{s-1+\varepsilon}\left(\mathbb{R}^{d}\right) \times H^{s}\left(\mathbb{R}^{d}\right), s>\frac{d}{2}, 0<\varepsilon<1$ in \cite{FEFF1}. Chemin-McCormick-Robinson-Rodrigo \cite{FEFF2} presented the local existence of weak solutions to \eqref{1.1} in $\mathbb{R}^{d}, d=2,3$ with the initial data $\left(u_{0}, b_{0}\right) \in B_{2,1}^{\frac{d}{2}-1}\left(\mathbb{R}^{d}\right) \times B_{2,1}^{\frac{d}{2}}\left(\mathbb{R}^{d}\right)$ and also proved the corresponding solution
is unique in $3 \mathrm{D}$ case. Wan in \cite{WAN} resolved the uniqueness of the solution in the 2D case by using mixed space-time Besov spaces. Recently,  Li-Tan-Yin \cite{LTY} made an important progress by reducing the functional setting to homogeneous Besov space $\left(u_{0}, b_{0}\right) \in \dot{B}_{p, 1}^{d-1}\left(\mathbb{R}^{d}\right) \times \dot{B}_{p, 1}^{\frac{d}{p}}\left(\mathbb{R}^{d}\right)$ where $p \in[1,2 d]$.

System \eqref{1.1} can be regarded as a coupling of a parabolic system with a hyperbolic one. For smooth solution, it (formally) obeys the following basic energy law:
\begin{equation}\label{1.2}
\frac{1}{2}\int_{\Omega}\left(|u|^{2}+|b|^{2}\right)(t) d x+\mu\int_{0}^{t} \int_{\Omega}|\nabla u|^{2} d x d s=\frac{1}{2} \int_{\Omega}\left(\left|u_{0}\right|^{2}+\left|b_{0}\right|^{2}\right) d x.
\end{equation}
However, for weak solutions with less regularity, \eqref{1.2} may fail. That is, a weak solution $(u,b)$ to \eqref{1.1} satisfies the following
energy inequality:
\begin{equation}\label{EI}
	\frac{1}{2}\int_{\Omega}\left(|u|^{2}+|b|^{2}\right)(t) d x+\mu\int_{0}^{t} \int_{\Omega}|\nabla u|^{2} d x d s\leq\frac{1}{2} \int_{\Omega}\left(\left|u_{0}\right|^{2}+\left|b_{0}\right|^{2}\right) d x.
\end{equation}
From the physical point of view we would expect that a weak solution $(u, b)$ satisfies
the energy equality \eqref{1.2}.
It is natural to ask the following question: Which regularity of weak solutions for the validity of the energy
equality \eqref{1.2}?

In \eqref{1.1}, if $\mu=0$, then the equations reduce to the ideal MHD equations. The famous Onsager conjecture for the Euler and ideal MHD equations predicts the
threshold regularity for energy conservation. In this direction,  Onsager \cite{ON} considers periodic 3-dimensional weak solutions of the incompressible Euler equations, where the velocity $u$ satisfies the uniform H\"older condition
$$
|u(x, t)-u(x', t)|\leq C|x-x'|^\alpha,
$$
for constants $C$ and $\alpha$ independent of $x, x'$ and $t$.\\
 $(a)$ If $ \alpha > \frac{1}{3}$
, then the total kinetic energy $E(t) = \frac{1}{2}\int |u(x, t)|^2 dx$ is a constant.\\
 $(b)$ For any $ \alpha < \frac{1}{3}$
, there are $u$ for which it is not a constant.\\
For $\alpha > \frac{1}{3}$, Eyink \cite{EGL} and Constantin-E-Titi \cite{CPT} proved the energy conservation. Recently,  Isett \cite{IS} gave the complete answer to part $(b)$ of the conjecture.
Cheskidov-Constantin-Friedlander-Shvydkoy in \cite{CCFS} proved that weak solution in the following class
\begin{equation}\label{1.3}
	u\in L^{3}\left(0, T; B_{3, c(\mathbb{N})}^{1 / 3}\right) \cap C_{w}\left(0, T; L^{2}\right)
\end{equation}
conserves energy. Bardos-Titi \cite{BT} considered the Onsager's conjecture for the incompressible
Euler equations in bounded domains and showed that, in a bounded domain $\Omega \subset \mathbb{R}^{d}$, with $\partial \Omega \in C^{2}$, any weak solution $(u(x, t), p(x, t))$, of the Euler equations of ideal incompressible fluid in $\Omega \times(0, T) \subset \mathbb{R}^{d} \times \mathbb{R}_{t}$, with the impermeability boundary condition $u \cdot \vec{n}=0$ on $\partial \Omega \times(0, T)$, is of constant energy on the interval $(0, T)$, provided the velocity field $u \in L^{3}\left((0, T) ; C^{0, \alpha}(\bar{\Omega})\right)$, with $\alpha>\frac{1}{3}$. Concerning the energy conservation of ideal MHD equations, by directly applying the methods developed in the study of Onsager's conjecture, several sufficient conditions have been obtained in \cite{CKS,KL,YUX}. More precisely, Caflisch-Klapper-Steele \cite{CKS} shown  that energy is conserved if
\begin{equation}\label{1.4}
	u\in C\left(0,T;B^{\alpha_1}_{3,\infty}\left(\mathbb{R}^3\right)\right);\,\, b\in C\left(0,T;B^{\alpha_2}_{3,\infty}\left(\mathbb{R}^3\right)\right)
\end{equation}
with $\alpha_1>\frac{1}{3}$ and $\alpha_1+2\alpha_2>1$. Later, sharper sufficient condition have been established in \cite{KL}:
\begin{equation}\label{1.5}
	u\in L^{3}\left(0,T;B^{\alpha_1}_{3,c(N)}\left(\mathbb{R}^3\right)\right);\,\, b\in L^{3}\left(0,T;B^{\alpha_2}_{3,c(N)}\left(\mathbb{R}^3\right)\right)
\end{equation}
with $\alpha_1\geq\frac{1}{3}$ and $\alpha_1+2\alpha_2\geq1$. Later, Yu \cite{YUX} improved those conditions and obtained new sufficient conditions for the conservation of energy.

%Besides their physical applications, the MHD equations \eqref{1.1}
%are also mathematically significant from the theoretical point of view. It is well known that the global existence of weak solutions, local
%existence, and uniqueness of smooth solutions to equations \eqref{1.1} were established in \cite{DL,TEMAM}.
%  a global weak solution to $(1.1)$ for initial data with
%finite energy, that is,
%\[u,b\in L^{\infty}\left(0,T;L^{2}(\mathbb{R}^3)\right)\cap L^{2}\left(0,T;H^{1}\left(\mathbb{R}^3\right)\right)\quad \text{for any}\quad T>0.\]
Notice that in the absence of magnetic field (i.e., $b=0$) and $\mu>0$, equations \eqref{1.1} reduce to the incompressible Navier-Stokes (NS) equations. As in the incompressible NS equations, the question of the regularity and uniqueness of weak solutions remains one of the biggest open problems in mathematical fluid mechanics. Energy equality is clearly a prerequisite for regularity, and can be a first step in proving conditional regularity results \cite{SS}.  Prodi \cite{GP},  Serrin \cite{SJ} and Escauriaza-Seregin-Sverak \cite{ELSGA} first proved conditional regularity result
\begin{equation}\label{1.6}
u\in L^{q}\left(0,T;L^{p}\left(\Omega\right)\right) \quad \text{with} \quad\frac{2}{q}+\frac{3}{p}\leq1 \quad \text{and}\quad 3\leq p\leq\infty,
\end{equation}
naturally, energy is conserved under this condition.
Lions \cite{LIONS} proved  that energy equality holds for weak solutions $u\in L^4\left(0, T; L^4(\mathbb{T}^d)\right)$. A few years later, Shinbrot \cite{SHI} improved upon this result to
\begin{equation}\label{1.7}
u\in L^{q}\left(0,T;L^{p}(\mathbb{T}^d)\right),\,\,  \frac{2}{q}+\frac{2}{p}\leq 1\,\, \text{with}\,\, p\geq 4.
\end{equation}
Note that the energy equality in the limiting
case (i.e., $u\in L^{2}\left(0,T;L^{\infty}\right)$) was generalized to $u\in L^{2}\left(0,T; BMO\right)$ by Kozono-Taniuchi \cite{KT}. Recent work by Leslie-Shvydkoy \cite{LS} established the energy
equality under new $L^{q}\left(0,T;L^{p}\right)$ conditions using local energy estimates. In another paper, they \cite{LS1} proved that any solution to the 3-dimensional Navier-Stokes Equations in $\mathbb{R}^{3}$ which is Type I (i.e. $\left.\|u(t)\|_{L^{\infty}} \leq \frac{C}{\sqrt{T-t}}\right)$ in time must satisfy the energy equality at the first blowup time $T$. It is worth pointing out that the methods developed in \cite{LS1} relies on a expanding type iterative argument which can not deal with the case
$ L^{2,\infty}\left(0,T; BMO\right)$ and domain with boundary. The energy equality at the first blowup time $T$, Cheskidov-Luo \cite{CLX} showed that the energy equality is holding in the weak-in-time Onsager spaces, as a corollary, they deduced that $u \in L^{q, \infty}\left([0, T] ; B_{p, \infty}^{0}\right)$ for $\frac{1}{q}+\frac{1}{p}=\frac{1}{2}, p>4$ implies energy equality. Recently, Tan-Yin \cite{TAN}  proved the energy equality holds if $u\in L^{2,\infty}\left(0,T; BMO\right)$ and pointed out that some potential Type II blow-up are admissible in $L^{2,\infty}\left(0,T; BMO\right)$.
For Shinbrot's result, recently, Yu \cite{YUCHENG} gave a new proof which relies on a crucial lemma introduced by Lions. All these results deal with either $ \Omega= \mathbb{R}^d$ or $\Omega= \mathbb{T}^d$. However, due to the well-recognized dominant role of the boundary in the generation of turbulence, it seems very reasonable to investigate the energy conservation in bounded domains. Cheskidov-Friedlander-Shvydkoy \cite{CF} proved energy equality for $u\in L^3D\left(A^\frac{5}{12} \right)$ on a bounded domain, here $A$ denotes the Stokes
operator. Recently, Berselli and Chiodaroli \cite{BLC} established some new energy balance criteria involving the gradient of the
velocity. Yu \cite{YC} proved the validity of energy equality under the following assumption
$$
u \in L^{q}\left([0, T] ; L^{p}(\Omega)\right) \cap L^{s}\left(0, T ; B_{s}^{\alpha, \infty}(\Omega)\right)
$$
for $\frac{1}{q}+\frac{1}{p} \leq \frac{1}{2}, p \geq 4$ and $\frac{1}{2}+\frac{1}{s}<\alpha<1, s>2$. Very recently, Chen-Liang-Wang-Xu \cite{CLWX} showed that the Shinbrot's condition \eqref{1.7} together with $P \in L^{2}\left(0, T ; L^{2}(\partial \Omega)\right)$ guaranteed the energy equality. The additional assumptions $u \in L^{s}\left(0, T ; B_{s}^{\alpha, \infty}(\Omega)\right)$ in \cite{YC} and $P \in L^{2}\left(0, T ; L^{2}(\partial \Omega)\right)$ in \cite{CLWX} are used to deal with the boundary effects, furthermore, they must assume that $\Omega$ is an open, bounded domain with $C^2$ boundary $\partial \Omega$.

For the standard MHD equations,  Yong-Jiu \cite{JIU} proved the energy equality of weak solutions under the assumption velocity $u$ belongs to the critical space $L^\infty(0,T; L^3(\Omega))$. It is interesting that the velocity field seems to play a dominant role in energy equality, but, it is not clear whether this condition can be relaxed into
\eqref{1.7}. More recently, Kim in \cite{KIM} applying Yu's argument established the energy equality for the MHD equations which allow velocity $u$ and magnetic field $b$ to belong to $L^{q}\left(0,T;L^{p}\left(\mathbb{T}^3\right)\right)$, where $\frac{2}{q}+\frac{2}{p}\leq 1$ with $ p\geq 4$, moreover, they obtained the lower bounds for possible singular solutions. And then, for the bounded domain, Wang and Zuo \cite{WANGZUO} established two types of regularity
conditions to ensure the energy equality, but they need to impose some additional conditions on pressure $P$ near the boundary to tackle the boundary effect.

We emphasize that the existed results about the problem of energy conservation are considering either hyperbolic system or parabolic system, there are rarely results to investigate the system \eqref{1.1} which couples a parabolic system of $u$ with a hyperbolic system of $b$. The purpose of this work is to consider the energy conservation of weak solutions for the non-resistive MHD equations \eqref{1.1} with physical boundaries.
%Here, the boundary $\partial\Omega$ is only required to be Lipschitz but not $C^2$. We would like to emphasize that the system \eqref{1.1} is only partially parabolic, owing to the
%magnetic equation which is of hyperbolic type.
Compared with the Euler or the Navier-Stokes equations, system \eqref{1.1} has stronger nonlinearity due to the couplings of the fluid velocity $u$ and the magnetic field $b$, in addition, due to equations \eqref{1.1} without magnetic diffusion, the only prior estimate we can get on the magnetic field is
$b \in L^{\infty}\left(0, T ; L^{2}\left(\Omega\right)\right) $. In spite of the difficulties due to the lack of magnetic diffusion, and the boundary effects are considered, we still remark that the energy conservation criteria on the magnetic field $b$ are essentially the same as that on the velocity $u$ %Our results require the regularity of boundary $\partial\Omega$ is only Lipschitz which is the minimum requirement to make the boundary condition $b\cdot n$ sense,
and there are no boundary layer assumptions and additional restrictions on pressure $P$.
Our critical strategies are that we noticing the important properties of weak solutions to the nonstationary Stokes system and the separate mollification of weak solutions from the boundary effect by considering a non-standard local energy equality and transform the boundary effects into the estimates of the gradient of cut-off functions. And then, we obtain the global energy equality by exploring the special nonlinear properties of equations \eqref{1.1} and taking suitable cut-off functions.

Now, we recall the definition of a Leray-Hopf weak solution.
\begin{definition}\label{def1}
Let $(u_0, b_0)\in L^2(\Omega)$ with $\nabla\cdot u_0=\nabla\cdot b_0 = 0$, $T > 0$. The function $(u, b)$ defined in $[0, T ] \times\Omega$ is said to be a Leray-Hopf weak solution to \eqref{1.1} if\\
1. $u \in L^{\infty}\left(0, T ; L^{2}\left(\Omega\right)\right) \cap L^{2}\left(0, T ; H^{1}\left(\Omega\right)\right), b \in L^{\infty}\left(0, T ; L^{2}\left(\Omega\right)\right)$; \\
2. for any $t \in[0, T], \Phi, \Psi \in C_{c}^{\infty}(\Omega \times[0, T])$, we have
$$
\begin{aligned}
	&\int_{\Omega} u(x, t) \cdot \Phi(x, t) d x-\int_{\Omega} u(x, 0) \cdot \Phi(x, 0) d x-\int_{0}^{t} \int_{\Omega} u(x, s) \cdot \partial_{s} \Phi(x, s) d x d s \\
	=&\int_{0}^{t} \int_{\Omega}(u \otimes u: \nabla \Phi-b \otimes b: \nabla \Phi+P d i v \Phi) d x d s-\mu \int_{0}^{T} \int_{\Omega} \nabla u: \nabla \Phi d x d s
\end{aligned}
$$
and
$$
\begin{aligned}
	&\int_{\Omega} b(x, t) \cdot \Psi(x, t) d x-\int_{\Omega} b(x, 0) \cdot \Psi(x, 0) d x-\int_{0}^{t} \int_{\Omega} b(x, s) \cdot \partial_{s} \Psi(x, s) d x d s \\
	&=\int_{0}^{t} \int_{\Omega}(b \otimes u: \nabla \Psi-u \otimes b: \nabla \Psi) d x d s
\end{aligned}
$$
3. for any $\varphi \in C_{c}^{\infty}(\Omega)$, it holds that
$$
\int_{\Omega} u \cdot \nabla \varphi d x=\int_{\Omega} b \cdot \nabla \varphi d x=0
$$
a.e. $t \in(0, T)$;\\
4. $(u,b)$ satisfies the energy inequality
$$	\frac{1}{2}\int_{\Omega}\left(|u|^{2}+|b|^{2}\right)(t) d x+\mu\int_{0}^{t} \int_{\Omega}|\nabla u|^{2} d x d s\leq\frac{1}{2} \int_{\Omega}\left(\left|u_{0}\right|^{2}+\left|b_{0}\right|^{2}\right) d x,$$
for all $t \in[0, T]$. Furthermore, it holds that $$\lim_{t\to 0^+}\left(\|u(\cdot,t)\|^2_{L^2}+\|b(\cdot,t)\|^2_{L^2}\right)=\|u_0\|^2_{L^2}+\|b_0\|^2_{L^2}.$$
\end{definition}

Now we state our main result.
\begin{theorem}\label{th1}Let $\Omega$ be a bounded Lipschitz domain in $\mathbb{R}^{d}$ with $d \geq 2$ and $(u, b)$ be a Leray-Hopf weak solution to the initial-boundary value problem \eqref{1.1} in the sense of Definition \ref{def1}. Assume that for any $0<\tau\leq T$
\begin{equation}\label{1.8}
	u\in L^{q}\left(\tau,T;L^{p}\left(\Omega\right)\right) \quad \text{with} \quad\frac{2}{q}+\frac{2}{p}\leq1, \quad  p\geq 4
\end{equation}
and
\begin{equation}\label{1.9}
	b\in L^{r}\left(\tau,T;L^{s}\left(\Omega\right)\right) \quad \text{with} \quad\frac{2}{r}+\frac{2}{s}\leq1, \quad  s\geq 4,
\end{equation}
then the energy equality \eqref{1.2} holds for all $t \in  [0, T ]$.
\end{theorem}
\begin{remark}
Our energy balance criteria also holds true for classical MHD equations without any difficulties. In \cite{WANGZUO}, they proved the following energy balance criteria:
$$u, b\in L^{4}((0, T) \times \Omega)
\,\, and \,\,
P \in L^{2}\left((0, T) \times \Omega_{\delta}\right),$$
where $\Omega_{\delta}=\{x \in \Omega \mid \operatorname{dist}(x, \partial \Omega)<\delta\}$. The novelty of this paper is that we consider the energy conservation for the non-resistive MHD equations, and we remark that our result needn't impose any condition on the pressure. It is also worth pointing out that our method also adapts to the case $\Omega=\mathbb{R}^d\,(d=2, 3, 4)$, and by using the geometric structure of $\mathbb{R}^d$, we can improve the Shinbrot's criteria into $u\in L^q_{loc}(0,T;L^p_{loc}(\mathbb{R}^d))$ with $\frac{1}{q}+\frac{1}{p}=\frac{1}{2}$, $4\leq p$ and
$b\in L^r_{loc}\left(0,T;L^s_{loc}\left(\mathbb{R}^d\right)\right)\cap L^{\frac{4d+8}{d+4}}\left(0,T;L^{\frac{4d+8}{d+4}}\left(\mathbb{R}^{d}\right)\right)$ with $\frac{1}{r}+\frac{1}{s}=\frac{1}{2}$, $4\leq s$.
\end{remark}
\begin{remark}
We will note that $$u \in L^{\infty}\left(0, T ; L^{2}(\Omega)\right) \cap L^{q}\left(\tau, T ; L^{p}(\Omega)\right)\text{ for any } \frac{1}{q}+\frac{1}{p} \leq \frac{1}{2},\,\, p \geq 4,$$
and
$$b \in L^{\infty}\left(0, T ; L^{2}(\Omega)\right) \cap L^{r}\left(\tau, T ; L^{s}(\Omega)\right) \text{ for any } \frac{1}{r}+\frac{1}{s} \leq \frac{1}{2},\,\,s \geq 4,$$
then one can deduce that
$$
\|u\|_{L^{4}\left(\tau, T ; L^{4}(\Omega)\right)} \leq C\|u\|_{L^{\infty}\left(\tau, T ; L^{2}(\Omega)\right)}^{a}\|u\|_{L^{q}\left(\tau, T ; L^{p}(\Omega)\right)}^{1-a}
$$
and
$$
\|b\|_{L^{4}\left(\tau, T ; L^{4}(\Omega)\right)} \leq C\|b\|_{L^{\infty}\left(\tau, T ; L^{2}(\Omega)\right)}^{a}\|b\|_{L^{r}\left(\tau, T ; L^{s}(\Omega)\right)}^{1-a}
$$
for some $0<a<1$. Thus, we can use the facts that $u, b$ are bounded in $L^{4}\left(\tau, T ; L^{4}(\Omega)\right)$ in our proof.
\end{remark}
\section{Preliminaries}
To investigate the boundary effect, we will recall some useful notations and Lemmas.\\
Assume $\Omega$ be a bounded domain with boundary $\partial \Omega$ in $\mathbb{R}^{d}$ with $d \geq 2$, for any $x\in\Omega$, $d(x)=\inf_{y\in\partial\Omega}|x-y|$, we define $\Omega^*_{\delta}$ as follows
$$
\Omega^*_{\delta}=\{x\in\Omega|d(x)>\delta\}.
$$
It is clear that $\Omega^*_{\delta}$ is still a domain if $\delta$ is small enough. We use $L^{q}\left(0, T ; L^{p}(\Omega)\right)$ to denote the space of measurable functions with the following norm
$$
\|f\|_{L^{q}\left(0, T ; L^{p}(\Omega)\right)}=\left\{\begin{array}{l}
	\left(\int_{0}^{T}\left(\int_{\Omega}|f(t, x)|^{p} d x\right)^{\frac{q}{p}} d t\right)^{\frac{1}{q}}, 1 \leq q<\infty \\
	\operatorname{ess} \sup _{t \in[0, T]}\|f(t, \cdot)\|_{L^{p}(\Omega)}, q=\infty.
\end{array}\right.
$$
 If $u\in L^q(\tau,T;L^p(\Omega))$ for any $0<\tau<T$, we say $u\in L^q_{loc}(0,T;L^p(\Omega))$.

Let $\eta: \mathbb{R}^{d} \rightarrow \mathbb{R}$ be a standard mollifier, i.e. $\eta(x)=C \mathrm{e}^{\frac{1}{|x|^{2}-1}}$ for $|x|<1$ and $\eta(x)=0$ for $|x| \geqslant 1$, where  constant $C>0$ selected such that $\int_{\mathbb{R}^{d}} \eta(x) \mathrm{d} x=1$. For any $\varepsilon>0$, we define the rescaled mollifier $\eta_{\varepsilon}(x)=\varepsilon^{-d} \eta\left(\frac{x}{\varepsilon}\right) .$ For any function $f \in L_{\text {loc }}^{1}(\Omega)$, its mollified version is defined as
$$f^{\varepsilon}(x)=\left(f * \eta_{\varepsilon}\right)(x)=\int_{\Omega} \eta_{\varepsilon}(x-y) f(y) \mathrm{d} y .$$
If $f \in  W^{1, p}(\Omega)$, the following local approximation is well known
$$
f^{\varepsilon}(x) \rightarrow f \quad \text { in } \quad  W_{l o c}^{1, p}(\Omega) \quad \forall p \in[1, \infty).
$$
The crucial ingredients to prove Theorem \ref{th1} are the following important lemmas. In the first lemma, we shall construct some cut-off functions to separate the mollification of weak solution from the boundary effect.
\begin{lemma}\label{CUT-OFF}
Let $\Omega$ be a bounded domain with boundary $\partial \Omega$. For any $\delta>0$ small enough, there exist some cut-off functions $\varphi(x)$ which satisfy $\varphi(x)=1$ for $x \in \Omega^*_{3 \delta}$ and $\varphi(x)=0$ for $x \in \Omega \backslash \Omega^*_{\delta}$. Furthermore, it holds that
$$
|\nabla \varphi(x)| \leq \frac{C}{\delta}.
$$
\textbf{Proof.} We first take $\chi(x)$ is the characteristic function of $\Omega^*_{2 \delta}$. We claim that for any $x \in \Omega^*_{2 \delta}$, it yields that
$$
\inf _{y \in \partial \Omega^*_{\delta}}|x-y| \geq \delta.
$$
We thus define $\varphi(x)= (\chi*\eta_{\delta})(x)$, it is obvious that $\varphi(x)=1$ for $x \in \Omega^*_{3 \delta}$ and $\varphi(x)=0$ for $x \in \Omega \backslash \Omega^*_{\delta}$
and $|\nabla \varphi(x)| \leq \frac{C}{\delta}$. We now prove our claim, for any $x \in \Omega^*_{2 \delta}$, assume $y_{0} \in \partial \Omega^*_{\delta}$ satisfy $\left|x-y_{0}\right|=$ $\inf _{y \in \partial \Omega^*_{\delta}}|x-y| .$ Similarly, we assume $z_{0} \in \partial \Omega$ satisfy $\left|y_{0}-z_{0}\right|=\inf _{z \in \partial \Omega}\left|y_{0}-z\right| .$ By the definitions
of $\Omega^*_{\delta}$ and $\Omega^*_{2 \delta}$, we have
$$
\left|x-z_{0}\right| \geq 2 \delta,\left|y_{0}-z_{0}\right|=\delta.
$$
This means that $\inf _{y \in \partial \Omega^*_{\delta}}|x-y|=\left|x-y_{0}\right| \geq\left|x-z_{0}\right|-\left|y_{0}-z_{0}\right| \geq \delta$.
\end{lemma}
The next lemma is the Lions commutator estimate.
\begin{lemma}\label{le2.1} \cite{LIONSJ}
Let $\partial$ be a partial derivative in one direction. Let $f, \partial f \in L^{p}\left(\mathbb{R}^{+} \times \Omega\right)$, $g \in L^{q}\left(\mathbb{R}^{+} \times \Omega\right)$ with $1 \leq p, q \leq \infty$, and $\frac{1}{p}+\frac{1}{q} \leq 1 .$ Then, we have
$$
\|\partial(f g)* \eta_{\varepsilon}-\partial\left(f (g* \eta_{\varepsilon})\right)\|_{L_{loc}^{r}\left(\mathbb{R}^{+} \times \Omega\right)} \leq C\|\partial f\|_{L^{p}\left(\mathbb{R}^{+} \times \Omega\right)}\|g\|_{L^{q}\left(\mathbb{R}^{+} \times \Omega\right)}
$$
for some constant $C>0$ independent of $\varepsilon, f$ and $g$, and with $\frac{1}{r}=\frac{1}{p}+\frac{1}{q} .$ In addition,
$$
\partial(f g)* \eta_{\varepsilon}-\partial\left(f (g* \eta_{\varepsilon})\right) \rightarrow 0 \quad \text { in } L_{loc}^{r}\left(\mathbb{R}^{+} \times \Omega\right)
$$
as $\varepsilon \rightarrow 0$, if $r<\infty$.
\end{lemma}

\begin{lemma}\label{le2.3}\cite{SOHR}
Let $\Omega \subseteq \mathbb{R}^{d}, d \geq 2$, be any domain, let $0<T \leq \infty$, $u_{0} \in$ $L_{\sigma}^{2}(\Omega), f=f_{0}+\operatorname{div} F$ with
$$
f_{0} \in L^{1}\left(0, T ; L^{2}(\Omega)\right), \quad F \in L^{2}\left(0, T ; L^{2}(\Omega)\right),
$$
and let
$$
u \in L_{l o c}^{1}\left([0, T) ; W_{0, \sigma}^{1,2}(\Omega)\right)
$$
be a weak solution of the Stokes system
$$
u_{t}-\mu \Delta u+\nabla P=f, \quad \operatorname{div} u=0,\left.\quad u\right|_{\partial \Omega}=0, \quad u(0)=u_{0}
$$
with data $f, u_{0}$.
Then $u$ has the following properties:

a) $u \in L^{\infty}\left(0, T ; L_{\sigma}^{2}(\Omega)\right), \quad \nabla u \in L^{2}\left(0, T ; L^{2}(\Omega)\right) .$

b) $u:[0, T) \rightarrow L_{\sigma}^{2}(\Omega)$ is strongly continuous, after a redefinition on a null set of $[0, T), u(0)=u_{0}$, and the energy equality
$$
\begin{aligned}
	\frac{1}{2}\|u(t)\|_{2}^{2}+\mu \int_{0}^{t}\|\nabla u\|_{2}^{2} d \tau=& \frac{1}{2}\left\|u_{0}\right\|_{2}^{2}+\int_{0}^{t}\langle f_{0}, u\rangle_{\Omega} d \tau -\int_{0}^{t}\langle F, \nabla u\rangle_{\Omega} d \tau.
\end{aligned}
$$	
\end{lemma}
\section{Proof of Theorem \ref{th1}}
The main object of this section is to prove Theorem \ref{th1}. For the sake of simplicity, we will proceed as if the solution is differentiable in time. The extra arguments needed to mollify in time are straightforward. First,  we rewrite the first equation of the system \eqref{1.1} as
$$\partial_t u-\mu\Delta u+\nabla P=\operatorname{div}\left((b\otimes b)-(u\otimes u)\right),$$
since
$$\left((b\otimes b)-(u\otimes u)\right) \in L_{l o c}^{2}\left((0, T) ; L^{2}(\Omega)\right),$$
it follows from Lemma \ref{le2.3} that the weak solution $u$ sastisfies
\begin{equation}
	\begin{aligned}\label{3.1}
		\frac{1}{2}\|u(t)\|_{L^2}^{2}+\mu\int_{0}^{t}\|\nabla u\|_{L^2}^{2} d \tau= \frac{1}{2}\left\|u_{0}\right\|_{L^2}^{2} -\int_{0}^{t}\langle (b\otimes b)-(u\otimes u), \nabla u\rangle_{\Omega} d \tau.
	\end{aligned}
\end{equation}
Next, we mollify the second equation of system \eqref{1.1},  then using $b^{\varepsilon}(t, x) \varphi(x)$ to test the result equation, integrating by parts over $[\tau, t] \times \Omega$ with $0<\tau \leq t \leq T$, one has
%$$
%	\begin{aligned}
%&\frac{1}{2} \int_{\Omega}\left|u^{\varepsilon}(t, x)\right|^{2} \varphi(x) d x-\frac{1}{2} \int_{\Omega}\left|u^{\varepsilon}(\tau, x)\right|^{2} \varphi(x) d x+\mu\int_{\tau}^{t} \int_{\Omega}\left|\nabla u^{\varepsilon}(s, x)\right|^{2} \varphi(x) d x d s\\
%=&\int_{\tau}^{t} \int_{\Omega}\left\{(u \otimes u)^{\varepsilon}: \nabla\left(u^{\varepsilon} \varphi\right)-(b \otimes b)^{\varepsilon}: \nabla\left(u^{\varepsilon} \varphi\right)+P^{\varepsilon} \nabla \cdot\left(u^{\varepsilon} \varphi\right)-\nabla u^{\varepsilon}:\left(u^{\varepsilon} \otimes \nabla \varphi\right)\right\} d x d s
%	\end{aligned}
%$$
%and
\begin{equation}
	\begin{aligned}\label{3.2}
		&\frac{1}{2} \int_{\Omega}\left|b^{\varepsilon}(t, x)\right|^{2} \varphi(x) d x-\frac{1}{2} \int_{\Omega}\left|b^{\varepsilon}(\tau, x)\right|^{2} \varphi(x) d x\\
		=&-\int_{\tau}^{t} \int_{\Omega}(u\cdot\nabla b)^{\varepsilon}b^{\varepsilon}\varphi(x) d x d s+\int_{\tau}^{t} \int_{\Omega}(b\cdot\nabla u)^{\varepsilon}b^{\varepsilon}\varphi(x) d x d s\\
		=&\mathcal{E}^{\varepsilon}(t),
	\end{aligned}
\end{equation}
where $\varphi(x)$ is a cut-off function constructed in Lemma \ref{CUT-OFF} and which equals to one on $\Omega^*_{3 \delta}$ and vanishes out of $\Omega^*_{\delta}$.
% Adding the above two equations together, we obtain
%\begin{equation}\label{3.9}
%\begin{aligned}
%	& \frac{1}{2}\int_{\Omega}\left(\left|u^{\varepsilon}(t, x)\right|^{2}\varphi(x)+\left|b^{\varepsilon}(t, x)\right|^{2}\varphi(x)\right) d x-\frac{1}{2}\int_{\Omega}\left(\left|u^{\varepsilon}(\tau, x)\right|^{2}\varphi(x)+\left|b^{\varepsilon}(\tau, x)\right|^{2}\varphi(x)\right) dx\\
%	&+\mu \int_{\tau}^{t}  \int_{\Omega}\left|\nabla u^{\varepsilon}(s,x)\right|^{2} \varphi(x)d x ds\\
%	=& \int_{\tau}^{t} \int_{\Omega}\left\{(u \otimes u)^{\varepsilon}: \nabla\left(u^{\varepsilon} \varphi\right)+P^{\varepsilon} \nabla \cdot\left(u^{\varepsilon} \varphi\right)-\nabla u^{\varepsilon}:\left(u^{\varepsilon} \otimes \nabla \varphi\right)\right\}dxds\\
%	&-\int_{\tau}^{t} \int_{\Omega}\left\{(b \otimes b)^{\varepsilon}: \nabla\left(u^{\varepsilon} \varphi\right)-(b\cdot\nabla u)^{\varepsilon} b^{\varepsilon}\varphi+(u\cdot\nabla b)^{\varepsilon} b^{\varepsilon}\varphi\right\} d x d s \\
%	=: & \mathcal{E}^\varepsilon(t).
%\end{aligned}
%\end{equation}
Note that
%$$
%\begin{aligned}
%(u \otimes u)^{\varepsilon}_{ij} &=[(u \otimes u)_{ij}^{\varepsilon}-(u \otimes u^{\varepsilon})_{ij}]+[(u \otimes u^{\varepsilon})_{ij}-(u^{\varepsilon} \otimes u^{\varepsilon})_{ij}]+(u^{\varepsilon} \otimes u^{\varepsilon})_{ij} \\
% &=[(u_iu_j)^{\varepsilon}-u_iu^{\varepsilon}_j]+[u_i  u^{\varepsilon}_j-u^{\varepsilon}_i u^{\varepsilon}_j]+u^{\varepsilon}_i u^{\varepsilon}_j,
%\end{aligned}
%$$
%$$
%\begin{aligned}
%	(b \otimes b)^{\varepsilon}_{ij} &=[(b \otimes b)_{ij}^{\varepsilon}-(b \otimes b^{\varepsilon})_{ij}]+[(b \otimes b^{\varepsilon})_{ij}-(b^{\varepsilon} \otimes b^{\varepsilon})_{ij}]+(b^{\varepsilon} \otimes b^{\varepsilon})_{ij} \\
%	&=[(b_ib_j)^{\varepsilon}-b_ib^{\varepsilon}_j]+[b_i  b^{\varepsilon}_j-b^{\varepsilon}_i b^{\varepsilon}_j]+b^{\varepsilon}_i b^{\varepsilon}_j
%\end{aligned}
%$$
$$
\begin{aligned}
	(u\cdot \nabla b)^{\varepsilon}_i=(u_j\partial_j b_i )^{\varepsilon} &=[\partial_j(u_jb_i)^{\varepsilon}-\partial_j(u_j
	b^{\varepsilon}_i)]+u_j\partial_jb^{\varepsilon}_i,
\end{aligned}
$$
$$
\begin{aligned}
(b\cdot \nabla u )^{\varepsilon}_i&=[\partial_j(u_i b_j )^{\varepsilon}-\partial_j(u_i b^{\varepsilon}_j )]+[\partial_j(u_i b^{\varepsilon}_j )-\partial_j(u^{\varepsilon}_i
	b^{\varepsilon}_j)]+b^{\varepsilon}_j\partial_ju^{\varepsilon}_i.
\end{aligned}
$$
It is clear that
%\begin{equation}\label{3.10}
%	\begin{aligned}
%	\mathcal{E}^{\varepsilon}(t)
%		=&\int_{\tau}^{t} \int_{\Omega}\left\{[(u_iu_j)^{\varepsilon}-u_iu^{\varepsilon}_j]+[u_i  u^{\varepsilon}_j-u^{\varepsilon}_i u^{\varepsilon}_j]+u^{\varepsilon}_i u^{\varepsilon}_j\right\}\partial_j(u^{\varepsilon}_i\varphi) dxds\\
%		&+ \int_{\tau}^{t} \int_{\Omega}\left\{P^{\varepsilon} \nabla \cdot\left(u^{\varepsilon} \varphi\right)-\nabla u^{\varepsilon}:\left(u^{\varepsilon} \otimes \nabla \varphi\right)\right\}dxds\\
%		&-\int_{\tau}^{t} \int_{\Omega}\left\{[(b_ib_j)^{\varepsilon}-b_ib^{\varepsilon}_j]+[b_i  b^{\varepsilon}_j-b^{\varepsilon}_i b^{\varepsilon}_j]+b^{\varepsilon}_i b^{\varepsilon}_j\right\}\partial_j(u^{\varepsilon}_i\varphi) dxds\\
%		&+\int_{\tau}^{t} \int_{\Omega}\left\{[\partial_j(u_i b_j )^{\varepsilon}-\partial_j(u_ib^{\varepsilon}_j )]+[\partial_j(u_i b^{\varepsilon}_j )-\partial_j(u^{\varepsilon}_i
%		b^{\varepsilon}_j)]+b^{\varepsilon}_j\partial_ju^{\varepsilon}_i\right\}b^{\varepsilon}_i
%		\varphi dxds\\
%		&-\int_{\tau}^{t} \int_{\Omega}\left\{[\partial_j(u_jb_i)^{\varepsilon}-\partial_j(u_j
%		b^{\varepsilon}_i)]+u_j\partial_jb^{\varepsilon}_i\right\}b^{\varepsilon}_i\varphi dxds.
%	\end{aligned}
%\end{equation}
\begin{equation}\label{3.3}
	\begin{aligned}
		\mathcal{E}^{\varepsilon}(t)
		=&\int_{\tau}^{t} \int_{\Omega}\left\{[\partial_j(u_i b_j )^{\varepsilon}-\partial_j(u_ib^{\varepsilon}_j )]+[\partial_j(u_i b^{\varepsilon}_j )-\partial_j(u^{\varepsilon}_i
		b^{\varepsilon}_j)]+b^{\varepsilon}_j\partial_ju^{\varepsilon}_i\right\}b^{\varepsilon}_i
		\varphi dxds\\
		&-\int_{\tau}^{t} \int_{\Omega}\left\{[\partial_j(u_jb_i)^{\varepsilon}-\partial_j(u_j
		b^{\varepsilon}_i)]+u_j\partial_jb^{\varepsilon}_i\right\}b^{\varepsilon}_i\varphi dxds.
	\end{aligned}
\end{equation}
In addition, using integration by parts
%$$\int_{\tau}^{t} \int_{\Omega}u^{\varepsilon}_i u^{\varepsilon}_j\partial_j(u^{\varepsilon}_i\varphi) dxds=\frac{1}{2}\int_{\tau}^{t} \int_{\Omega}|u^{\varepsilon}_i|^2u^{\varepsilon}_j\partial_j\varphi dxds,$$
%$$\int_{\tau}^{t} \int_{\Omega}b^{\varepsilon}_i b^{\varepsilon}_j\partial_j(u^{\varepsilon}_i\varphi) dxds=-\int_{\tau}^{t} \int_{\Omega}b^{\varepsilon}_j\partial_jb^{\varepsilon}_iu^{\varepsilon}_i
%\varphi dxds,$$
$$\int_{\tau}^{t} \int_{\Omega}u_j \partial_jb^{\varepsilon}_ib^{\varepsilon}_i\varphi dxds=-\frac{1}{2}\int_{\tau}^{t} \int_{\Omega}u_j\partial_j\varphi|b^{\varepsilon}_i|^2
 dxds,$$
which implies
\begin{equation}\label{3.4}
	\begin{aligned}
		\mathcal{E}^{\varepsilon}(t)
		=&\int_{\tau}^{t} \int_{\Omega}\left\{[\partial_j(u_i b_j )^{\varepsilon}-\partial_j(u_ib^{\varepsilon}_j )]+[\partial_j(u_i b^{\varepsilon}_j )-\partial_j(u^{\varepsilon}_i
		b^{\varepsilon}_j)]+b^{\varepsilon}_j\partial_ju^{\varepsilon}_i\right\}b^{\varepsilon}_i
		\varphi dxds\\
		&+\frac{1}{2}\int_{
			\tau}^{t} \int_{\Omega}u_j\partial_j\varphi |b^{\varepsilon}_i|^2dxds-\int_{\tau}^{t} \int_{\Omega}[\partial_j(u_jb_i)^{\varepsilon}-\partial_j(u_j
		b^{\varepsilon}_i)]b^{\varepsilon}_i\varphi dxds.
	\end{aligned}
\end{equation}
Combining \eqref{3.2} with \eqref{3.4}, we find
\begin{equation}
	\begin{aligned}\label{3.5}
		&\frac{1}{2} \int_{\Omega}\left|b^{\varepsilon}(t, x)\right|^{2} \varphi(x) d x-\frac{1}{2} \int_{\Omega}\left|b^{\varepsilon}(\tau, x)\right|^{2} \varphi(x) d x\\
		=&\int_{\tau}^{t} \int_{\Omega}\left\{[\partial_j(u_i b_j )^{\varepsilon}-\partial_j(u_ib^{\varepsilon}_j )]+[\partial_j(u_i b^{\varepsilon}_j )-\partial_j(u^{\varepsilon}_i
		b^{\varepsilon}_j)]\right\}b^{\varepsilon}_i
		\varphi dxds\\
		&+\frac{1}{2}\int_{
			\tau}^{t} \int_{\Omega}u_j\partial_j\varphi |b^{\varepsilon}_i|^2dxds-\int_{\tau}^{t} \int_{\Omega}[\partial_j(u_jb_i)^{\varepsilon}-\partial_j(u_j
		b^{\varepsilon}_i)]b^{\varepsilon}_i\varphi dxds\\
		&+\int_{\tau}^{t} \int_{\Omega}b^{\varepsilon}_j\partial_ju^{\varepsilon}_ib^{\varepsilon}_i
		\varphi dx ds\\
		=&\mathcal{E}_1^{\varepsilon}+\mathcal{E}_2^{\varepsilon}+\mathcal{E}_3^{\varepsilon}+\mathcal{E}_4^{\varepsilon},
	\end{aligned}
\end{equation}
To bound $\mathcal{E}^{\varepsilon}_1$, by H\"older's inequality Lemma \ref{le2.1}, one easily deduces that
$$
\begin{aligned}
	\mathcal{E}^{\varepsilon}_1\leq&\left|\int_{\tau}^{t} \int_{\Omega}\left\{[\partial_j(u_i b_j )^{\varepsilon}-\partial_j(u_i b^{\varepsilon}_j )]+[\partial_j(u_i b^{\varepsilon}_j )-\partial_j(u^{\varepsilon}_i
	b^{\varepsilon}_j)]\right\}b^{\varepsilon}_i\varphi dxds\right| \\
	\leq &C\left\|\partial_j\left(u_{i} b_{j }\right)^{\varepsilon}-\partial_j(u_ib^{\varepsilon}_j )\right\|_{L^{\frac{4}{3}}\left(\tau, T ; L^{\frac{4}{3}}(\Omega^*_{\delta})\right)}\left\|(b^{\varepsilon}_i\varphi) \right\|_{L^{4}\left(\tau, T ; L^{4}(\Omega^*_{\delta})\right)} \\
	&+C\left\|\partial_ju_{i}-\partial_ju^{\varepsilon}_i \right\|_{L^{2}\left(\tau, T ; L^{2}(\Omega^*_{\delta})\right)}\left\|(|b^{\varepsilon}_i|^2\varphi) \right\|_{L^{2}\left(\tau, T ; L^{2}(\Omega^*_{\delta})\right)}\\
	\rightarrow & 0, \quad as\quad \epsilon\rightarrow 0.
\end{aligned}
$$
By using the same trick as $\mathcal{E}^{\varepsilon}_1$, we have
$$
\begin{aligned}
	\mathcal{E}^{\varepsilon}_3\leq&\left|-\int_{\tau}^{t} \int_{\Omega}[\partial_j(u_jb_i)^{\varepsilon}-\partial_j(u_j
	b^{\varepsilon}_i)]b^{\varepsilon}_i\varphi dxds\right| \\
	\leq &C\left\|\partial_j\left(u_{j} b_{i }\right)^{\varepsilon}-\partial_j(u_j b^{\varepsilon}_i )\right\|_{L^{\frac{4}{3}}\left(\tau, T ; L^{\frac{4}{3}}(\Omega^*_{\delta})\right)}\left\|(b^{\varepsilon}_i\varphi) \right\|_{L^{4}\left(\tau, T ; L^{4}(\Omega^*_{\delta})\right)} \\
	\rightarrow & 0, \quad as\quad \varepsilon\rightarrow 0.
\end{aligned}
$$
For the term $\mathcal{E}^{\varepsilon}_4$, we claim that
%对于\,$\mathcal{E}_4^{\varepsilon}$, 我们断言
$$
\mathcal{E}_4^{\varepsilon}=\int_{\tau}^{t} \int_{\Omega}b^{\varepsilon}_j\partial_ju^{\varepsilon}_ib^{\varepsilon}_i
\varphi dx ds\rightarrow \int_{\tau}^{t} \int_{\Omega}b_j\partial_ju_ib_i
\varphi dx ds,\,\, as \,\,\varepsilon \rightarrow 0.
$$
Indeed,
\begin{equation*}
	\begin{split}
		&\int_{\tau}^{t} \int_{\Omega}b^{\varepsilon}_j\partial_ju^{\varepsilon}_ib^{\varepsilon}_i
		\varphi -b_j\partial_ju_ib_i
		\varphi dxds\\
		=&\int_{\tau}^{t} \int_{\Omega}(b^{\varepsilon}_j-b_j)\partial_ju^{\varepsilon}_ib^{\varepsilon}_i
		\varphi+b_j(\partial_ju^{\varepsilon}_i-\partial_j u_i)b^{\varepsilon}_i
		\varphi+b_j\partial_j u_i(b^{\varepsilon}_i-b_i)
		\varphi dxds\\
		\leq&\left\|\partial_ju_{i}-\partial_ju^{\varepsilon}_i \right\|_{L^{2}\left(\tau, T ; L^{2}(\Omega^*_{\delta})\right)}\left\|b^{\varepsilon}_ib_j\varphi \right\|_{L^{2}\left(\tau, T ; L^{2}(\Omega^*_{\delta})\right)}\\
		&+\left\|b_j-b^{\varepsilon}_j \right\|_{L^{4}\left(\tau, T ; L^{4}(\Omega^*_{\delta})\right)}\left\|\partial_ju^{\varepsilon}_ib^{\varepsilon}_i\varphi \right\|_{L^{\frac{4}{3}}\left(\tau, T ; L^{\frac{4}{3}}(\Omega^*_{\delta})\right)}\\
		&+\left\|b_i-b^{\varepsilon}_i \right\|_{L^{4}\left(\tau, T ; L^{4}(\Omega^*_{\delta})\right)}\left\|\partial_ju_ib_j\varphi \right\|_{L^{\frac{4}{3}}\left(\tau, T ; L^{\frac{4}{3}}(\Omega^*_{\delta})\right)}\\
		&\rightarrow 0,\,\, as \,\,\varepsilon \rightarrow 0,
	\end{split}
\end{equation*}
Therefore, we first take $\varepsilon\rightarrow 0$ in \eqref{3.5} to see that
\begin{equation}
	\begin{aligned}\label{3.6}
		&\frac{1}{2} \int_{\Omega}\left|b(t, x)\right|^{2} \varphi(x) d x-\frac{1}{2} \int_{\Omega}\left|b(\tau, x)\right|^{2} \varphi(x) d x\\
		=&\frac{1}{2}\int_{
			\tau}^{t} \int_{\Omega}u_j\partial_j\varphi |b_i|^2dxds+\int_{\tau}^{t} \int_{\Omega}b_j\partial_ju_ib_i
		\varphi dx ds,
	\end{aligned}
\end{equation}
Using Lemma \ref{CUT-OFF} and H\"older's inequality, we have
$$
\begin{aligned}
	\frac{1}{2}\int_{
		\tau}^{t} \int_{\Omega}u_j\partial_j\varphi |b_i|^2dxds
	\leq	 \frac{C}{\delta}\|b\|_{L^{4}\left(\tau, T ; L^{4}\left(\Omega \backslash \Omega^*_{3 \delta}\right)\right)}^{2}\|u\|_{L^{2}\left(\tau, T ; L^{2}\left(\Omega \backslash \Omega^*_{3 \delta}\right)\right)}.
\end{aligned}
$$
Noticing that $u|_{\partial \Omega}=0$ and $\Omega$ is a bounded domain with a Lipschitz boundary, by using the Poincar\'e inequality, we find that
$$
\|u\|_{L^{2}\left(\tau, T ; L^{2}\left(\Omega \backslash \Omega^*_{3 \delta}\right)\right)} \leq C \delta\|\nabla u\|_{L^{2}\left(\tau, T ; L^{2}\left(\Omega \backslash \Omega^*_{3 \delta}\right)\right)}.
$$
Obviously,
$$
\frac{1}{2}\int_{
	\tau}^{t} \int_{\Omega}u_j\partial_j\varphi |b_i|^2dxds \rightarrow 0, \quad \text { as } \delta \rightarrow 0.
$$
Letting $\delta$ go to zero in \eqref{3.6}, we obtain
\begin{equation}
	\begin{aligned}\label{3.7}
		&\frac{1}{2}\int_{\Omega}\left|b(t, x)\right|^{2} d x-\frac{1}{2}\int_{\Omega}\left|b(\tau, x)\right|^{2} dx=\int_{\tau}^{t} \int_{\Omega}b_j\partial_ju_ib_i dx ds, \quad 0<\tau\leq t\leq T.
	\end{aligned}
\end{equation}
Combining \eqref{3.1} and \eqref{3.7} gives
\begin{equation}
	\begin{aligned}\label{3.8}
		&\frac{1}{2}\left(\|u(\cdot,t)\|^2_{L^2}+\|b(\cdot,t)\|^2_{L^2}\right)+\int_{0}^{t}\|\nabla u\|_{L^2}^{2} d s\\
		=& \frac{1}{2}\left(\left\|u_{0}\right\|_{L^2}^{2}+\left\|b(\tau)\right\|_{L^2}^{2} \right)-\int_{0}^{t}\langle(b\otimes b)-(u\otimes u), \nabla u\rangle_{\Omega} d s+\int_{\tau}^{t} \int_{\Omega}b_j\partial_ju_ib_i dx ds.
	\end{aligned}
\end{equation}
By passing to the limit as $\tau\rightarrow 0$ and using the facts that $$\lim_{t\to 0^+}\left(\|u(\cdot,t)\|^2_{L^2}+\|b(\cdot,t)\|^2_{L^2}\right)=\|u_0\|^2_{L^2}+\|b_0\|^2_{L^2},$$ we immediately get
\begin{equation}
	\begin{aligned}\label{I}
		&\frac{1}{2}\int_{\Omega}\left(|u|^{2}+|b|^{2}\right)(t) d x+\mu\int_{0}^{t} \int_{\Omega}|\nabla u|^{2} d x d s\\
		=&\frac{1}{2} \int_{\Omega}\left(\left|u_{0}\right|^{2}+\left|b_{0}\right|^{2}\right) d x+\int_{0}^{t} \int_{\Omega}(u\otimes u) \nabla u dx ds.
	\end{aligned}
\end{equation}
Next, we will prove that
$$\int_{0}^{t} \int_{\Omega}(u\otimes u) \nabla u dx ds\equiv0.$$
Now, $u(t) \in H^{1}_{0}(\Omega)$, for a.a. $t \in[0, T)$ and so, for any such fixed $t$, denoting by $\left\{\psi_{k}\right\}$ a sequence from $\left\{\psi_{k}\right\} \subset\mathcal{D}(\Omega)=\left\{\boldsymbol{\psi} \in C_{0}^{\infty}(\Omega): \operatorname{div} \psi=0 \text { in } \Omega\right\}$ converging to $u$ in $H^{1}$ we have
$$
\begin{aligned}
	\left|\langle u \cdot \nabla u, u\rangle-\langle u \cdot \nabla \psi_{k}, \psi_{k}\rangle\right| & \leq\left|\langle u \cdot \nabla u,\left(u-\psi_{k}\right)\rangle\right|+\left|\langle u \cdot \nabla\left(u-\psi_{k}\right), u\rangle\right| \\
	& \leq\|u\|_{L^4}\|\nabla u\|_{L^2}\left\|u-\psi_{k}\right\|_{L^4}+\|u\|_{L^4}^{2}\left\|\nabla\left(u-\psi_{k}\right)\right\|_{L^2}
\end{aligned}
$$
so we deduce
$$
\lim _{k \rightarrow \infty}\langle u \cdot \nabla \psi_{k}, \psi_{k}\rangle=\langle u \cdot \nabla u, u\rangle.
$$
However, since $u(t) \in H(\Omega)$ for a.a. $t$, we get
$$
\langle u \cdot \nabla u, u\rangle=\lim_{k\to\infty}\langle u \cdot \nabla \psi_{k}, \psi_{k}\rangle=\lim_{k\to\infty}\frac{1}{2}\langle u, \nabla\left(\psi_{k}\right)^{2}\rangle=0.
$$
Furthermore, $\langle u\cdot \nabla u, u\rangle\in L^1([0,T])$, hence
$$\int_{0}^{t} \int_{\Omega}(u\otimes u) \nabla u dx ds=0.$$
Thus, the proof of  Theorem \ref{th1} is finished.
\section*{Acknowledgments}

The authors are supported by the Construct Program of the Key Discipline in Hunan Province and NSFC Grant No. 11871209.

\end{document}